\documentclass[journal,twoside,web]{ieeecolor}
\usepackage{generic}
\usepackage{cite}
\usepackage{amsmath,amssymb,amsfonts}
\usepackage{algorithmic}
\usepackage{graphicx}
\usepackage{algorithm,algorithmic}
\usepackage{hyperref}
\hypersetup{hidelinks=true}
\usepackage{textcomp}
\usepackage{booktabs}

\newtheorem{lemma}{Lemma}
\newtheorem{assumption}{Assumption}
\newtheorem{remark}{Remark}
\newtheorem{theorem}{Theorem}
\newtheorem{example}{Example}
\newtheorem{corollary}{Corollary}

\def\BibTeX{{\rm B\kern-.05em{\sc i\kern-.025em b}\kern-.08em
T\kern-.1667em\lower.7ex\hbox{E}\kern-.125emX}}
\begin{document}

\title{On the Robustness of Derivative-free Methods for Linear Quadratic Regulator}

\author{{Weijian Li, Panagiotis Kounatidis, and Zhong-Ping Jiang, and Andreas A. Malikopoulos}
\thanks{This research was supported in part by NSF under Grants CNS-2401007, CMMI-2348381, IIS-2415478, and in part by MathWorks.}
\thanks{Weijian Li, Panagiotis Kounatidis and  Andreas A. Malikopoulos are with the Department of Civil and Environmental Engineering, Cornell University, Ithaca, NY, USA (e-mail: wl779@cornell.edu, pk586@cornell.edu, amaliko@cornell.edu).}
\thanks{Zhong-Ping Jiang is with the Control and Networks Laboratory, Department of Electrical and Computer Engineering, Tandon School of Engineering, New York University, Brooklyn,
NY 11201 USA (e-mail: zjiang@nyu.edu).}}

\maketitle

\begin{abstract}
Policy optimization has emerged as a central approach in reinforcement learning, with derivative-free methods receiving significant attention for linear quadratic regulator (LQR) problems involving unknown system dynamics. In this paper, we analyze the robustness of such methods for infinite-horizon LQR settings, where policy gradients are estimated using sampled cost evaluations. We consider the impact of perturbations—arising from sources such as function approximation errors and measurement noise—on the convergence behavior of these estimators. Our analysis demonstrates that, under sufficiently small perturbations, the resulting derivative-free algorithms converge to an arbitrarily small neighborhood of the optimal policy. We further derive explicit bounds on the allowable perturbation magnitude and establish the sample complexity required to ensure convergence guarantees in the perturbed setting.
\end{abstract}


\section{Introduction}
\label{sec:introduction}

Reinforcement learning has attracted considerable attention \cite{sutton1998reinforcement}, and has been applied to the control of unknown systems \cite{kiumarsi2017optimal, Malikopoulos2009a,Malikopoulos2009b,Malikopoulos2011, nakka2022multi}.
Recently, there has been renewed interest in a particular instance of such problems, namely the linear quadratic regulator (LQR) problem, which involves controlling an (unknown) linear dynamical system with a quadratic cost \cite{lewis2012optimal,  rizvi2019reinforcement}. The problem is theoretically tractable and serves as a tool for the feedback design in engineering systems \cite{jiang2017robust, lewis2012reinforcement, malikopoulos2018decentralized}.
Recent research \cite{Malikopoulos2022a,Malikopoulos2024} has highlighted the importance of reconciling learning and control for systems with unknown or partially known dynamics, especially in contexts where data arrive sequentially and system models evolve over time. Motivated by these ideas, we explore the robustness properties of derivative-free policy optimization methods for infinite-horizon LQR problems, where the gradient of the cost function is estimated solely through sampled cost values that may be subject to perturbations.

Interestingly, many successes in reinforcement learning hinge on the development of policy optimization methods such as policy gradient methods \cite{sutton1999policy}, actor-critic methods \cite{grondman2012survey}, trust-region methods \cite{schulman2015trust}, where the policy, mapping states or observations to actions, is parametrized and directly optimized.
These methods are easy to be implemented, and can handle high-dimensional and continuous state/action spaces \cite{zhang2019policy}.
In addition, some observations indicate that they typically achieve faster convergence than value-based algorithms \cite{o2016combining}.
In contrast to their empirical success, the theoretical understanding on these methods remains mysterious.
A recent line of research has been devoted to the properties and limitations of policy optimization methods for LQR problems \cite{dean2020sample, hu2023toward}.
To be specific, the control policy is parameterized as a linear function of the state, and then, the problem is solved by policy gradient methods under the framework of constrained optimization.
For instance, the standard gradient descent, natural policy gradient, and Newton-type algorithms were proposed in \cite{fazel2018global} with convergence guarantees.
The authors of \cite{fatkhullin2021optimizing} investigated properties of both static state and output feedback policies for an LQR problem, and designed a gradient descent algorithm. 
Exponential stability of the gradient flow was established in  \cite{mohammadi2021convergence}.
In \cite{zhao2023global}, primal-dual policy gradient algorithms were developed for an LQR problem with risk constraints.
When the dynamics of the LQR model are unknown, the policy gradient is not analytically available and must instead be approximated through finite-sample evaluations of the associated cost function.
Inspired by the idea, derivative-free methods have been developed \cite{fazel2018global, mohammadi2021convergence, zhao2023global}.
In \cite{malik2020derivative}, the sample complexity of derivative-free algorithms was established.
Following that, a derivative-free algorithm was developed in \cite{li2021distributed} for a decentralized LQR problem.

In practice, accurate policy gradient estimators are difficult or even impossible to derive due to different kinds of errors from function approximations, measurement noises, and external disturbances. As a consequence, it is significant to investigate the robustness of policy gradient methods \cite{jiang2020learning}.
In \cite{pang2021robust}, it was shown that Kleinman’s policy iteration algorithm, equivalent to the Newton gradient algorithm with a step size of $1/2$, can find a near-optimal policy even in the presence of estimation errors.
For perturbed gradient flows, its property of the small-disturbance input-to-state stability (ISS) was established in \cite{sontag2022remarks} with applications to an LQR problem.
Specifically, trajectories of the flows eventually enter and remain within a small neighborhood of the optimum if perturbations are small enough.
Then it was shown in \cite{cui2024small} that
perturbed gradient flows for general nonlinear programming problems are also small-disturbance ISS.
Based on the small-disturbance ISS characterization, a robust reinforcement learning framework was proposed in \cite{cui2024robust} for a risk-sensitive linear quadratic Gaussian control problem.
For derivative-free methods, it is critical to obtain the exact costs by sampling and experimentation since policy gradients are estimated by 
cost evaluations \cite{malik2020derivative, fazel2018global, tang2020distributed, mohammadi2021convergence}, but this is challenging.
For instance, under a given policy, the cost of an infinite-horizon LQR problem is generally approximated by a finite-horizon cost \cite{fazel2018global}.
In a decentralized scenario, the global cost was estimated in a distributed manner with errors \cite{li2021distributed}.
However, as far as we know, the effect of perturbations has not been addressed on derivative-free methods in existing literature.

Motivated by the above observations, we investigate the robustness of derivative-free methods for solving an infinite-horizon LQR problem. Our main contributions are at least twofold.
a) We consider perturbed derivative-free methods that are closely related to the algorithms in \cite{malik2020derivative, fazel2018global, li2021distributed, tang2020distributed},  but our analysis incorporates perturbations.
To the best of our knowledge, this is the first work to investigate the effect of perturbations on derivative-free methods.
b) We show that as long as the perturbations are small enough, the perturbed derivative-free methods can converge to any pre-specified neighborhood of the optimal solution by choosing appropriate stepsize and smoothing radius. We derive the explicit perturbation bound and analyze the sample complexity of the perturbed methods.
Furthermore, we provide the rollout length for derivative-free methods to solve the LQR problem.

This paper is organized as follows. Section \uppercase\expandafter{\romannumeral2} formulates the LQR problem and introduces the perturbed derivative-free methods. Section \uppercase\expandafter{\romannumeral3} focuses on the convergence analysis for the perturbed methods.
Illustrative examples are provided in Section  \uppercase\expandafter{\romannumeral4}, and concluding remarks are given in  Section  \uppercase\expandafter{\romannumeral5}.

\textbf{Notation:}
Let $\mathbb R^n$, $\mathbb R^{m \times n}$ and $I_n$ be the set of $n$-dimensional real column vectors, the set of $m$-by-$n$ dimensional real matrices, and the $n$-by-$n$ identity matrix.
Denote by $\mathbb S_+^{n}$ ($\mathbb S_{++}^{n}$) the set of $n$-by-$n$ positive (semi-)definite matrices.
Let $\langle \cdot, \cdot\rangle_F$ be the  Frobenius inner product of two matrices, $\Vert \cdot \Vert_2$ be the $l_2$-norm of a vector or matrix, and $\Vert \cdot \Vert_F$ be the Frobenius norm of a matrix.
Let $(\cdot)^\top$ be the transpose, and ${\rm Tr}(\cdot)$ be the trace of a matrix.
Define
$\mathcal S^{m \times n} := \{D \in \mathbb R^{m \times n} : \Vert D\Vert_F = 1\}.$
Denote by ${\rm Unif}(\mathcal{S}^{m \times n})$ the uniform distribution over $\mathcal S^{m \times n}$.
Let $\mathbb E[\cdot]$ be the expectation of a random variable.
The indicator function of a random event $S$ is denoted by $1_{S}$ such that $1_{S} = 1$ if the event $S$ occurs, and $1_{S} = 0$ otherwise.

\section{PROBLEM FORMULATION}

In this section, we introduce an infinite-horizon LQR problem and perturbed derivative-free methods.

\subsection{LQR Formulation}

Consider a discrete-time linear time-invariant (LTI) system described by 
\begin{equation}
\label{LTI:sys}
x_{t+1} = A x_t + Bu_t, ~x_0 \sim \mathcal D,
\end{equation}
where $x_t \in \mathbb R^n$ is the state, $u_t \in \mathbb R^m$ is the control input, $x_0 \in \mathbb R^n$ is the initial state from a distribution $\mathcal D$, and $A$ and $B$ are unknown constant matrices with appropriate dimensions.
The infinite-horizon LQR problem entails seeking a control policy $u := \{u_t\}_{t = 0}^\infty$ that minimizes the cost given by
\begin{equation}
\label{cost}
\mathcal J(u) := \mathbb E \Big[\sum_{t = 0}^\infty c_t(x_t, u_t) \Big],
\end{equation}
where $c_t(x_t, u_t) := x_t^\top Q x_t + u_t^\top R u_t$ with $Q \in \mathbb R^{n \times n}$ and $R \in \mathbb R^{m \times m}$.

We make the following standing assumptions, which are always valid throughout this paper.
\begin{assumption}
\label{ass:mat}
$(A, B)$ is stabilizable, $Q \in \mathbb S_{+}^n$, $R \in \mathbb S_{++}^m$ and $(A, Q^{\frac 12})$ is detectable.
\end{assumption}

\begin{assumption}
\label{ass:noise}
For a random variable $v \sim \mathcal D$, we have
\begin{equation*}
\mathbb E[v] = 0, ~
\mathbb E[vv^\top] = I_n, ~{\rm and}~
\Vert v\Vert_2^2 \le C_m~{\rm a.s.}
\end{equation*}
where $C_m > 0$ is a constant.
\end{assumption}

Referring to \cite{anderson2007optimal}, the LQR problem admits a unique optimal controller $u^*_t = - K^* x_t$, where $K^* \in \mathbb R^{m \times n}$ satisfies
$$K^* = (B^\top P^* B + R)^{-1} B^\top P^* A,$$
and $P^* \in \mathbb S^n_{++}$ is the unique  solution to the discrete algebraic Riccati equation (DARE) as
$$
P^* = A^\top P^* A + Q - 
A^\top P^* B (B^\top P^* B + R)^{-1} B^\top P^* A.
$$

Denote by $\mathcal K_{st}$ the set of admissible stabilizing control policies, i.e.,
\begin{equation}
\label{Kst}
\mathcal K_{st} = \{K \in \mathbb R^{m \times n} \mid \rho(A - BK) < 1\},
\end{equation}
where $\rho(\cdot)$ is the spectral radius.
We reformulate the LQR problem by a linearly parameterized policy as $u_t = -Kx_t$.
Given a policy $K \in \mathcal K_{st}$, we let $\mathcal J(K; x_0)$ be the cost from the initial state $x_0$, and let $\mathcal J(K)$ be the cost of the LQR problem.
Plugging the policy $K$ into (\ref{LTI:sys}) and (\ref{cost}), we derive
\begin{equation*}
\begin{aligned}
\mathcal J(K; x_0) &= \sum_{t = 0}^\infty
x_0^\top (A \!-\! BK)^{t, \top} (Q \!+\! K^\top R K)(A \!-\! BK)^t x_0 \\
&= x_0^\top  P_K x_0,
\end{aligned}
\end{equation*}
where $P_K \in \mathbb S^n_{++}$ is the solution to
$$
P_K = Q + K^\top R K + (A - BK)^\top P_K (A - BK).
$$

It is clear that 
\begin{equation*}
\label{def:JK}	
\mathcal J(K) = \mathbb E[\mathcal J(K; x_0)] = \mathbb E_{x_0 \sim \mathcal D}[x_0^\top P_K x_0].
\end{equation*}
Recalling $\mathbb E[x_0 x_0^\top] = I_n$ yields
$\mathcal J(K) = {\rm Tr}(P_K)$.
For any $K \in \mathcal K_{st}$, $\mathcal J(K) \ge \mathcal J(K^*)$ due to $P_K \succeq P^*$. Besides, since $\mathcal J(K) = \mathcal J(K^*)$ implies 
$P_K = P^*$ and $K = K^*$, $\mathcal J(K)$ has a unique minimum at $K^*$.
Thus, the optimal policy $K^*$ can be obtained by solving the following optimization problem
\begin{equation}
\label{reform}
\min_{K \in \mathcal K_{st}} ~ \mathcal J(K).
\end{equation}

\subsection{Perturbed Derivative-free Methods}

Throughout this paper, we consider a model-free setting. Specifically, the model parameters of (\ref{LTI:sys}) are unknown, but for a given policy, the cost is accessible.
In this case, a typical policy gradient method is
\begin{equation}
\label{alg:dev:free}
K_{s+1} = K_s - \eta G(K_s),
\end{equation}
where $G(K_s)$ is an estimator of the gradient $\nabla \mathcal J(K_s)$,
$\eta > 0$ is the stepsize, and the initial policy $K_0$ should lie in $\mathcal K_{st}$.
The gradient estimator $G(K_s)$ is computed by the costs at finitely many points.
Similar to \cite{malik2020derivative, li2021distributed, fazel2018global}, for a given initial state $x_0$ and a policy $K$, we consider the one-point gradient estimate given by
\begin{equation}
\label{one_point}
G_r^1(K, D; x_0) :=  \mathcal J(K + rD; x_0) \frac{d}{r} D,
\end{equation}
and its corresponding two-point analogue as
\begin{equation}
\label{two_point}
G_r^2(K, D; x_0) := [\mathcal J(K + rD; x_0) - \mathcal J(K - rD; x_0)]\frac {d}{2r} D,
\end{equation}
where $d = mn$, $r > 0$ is the smoothing radius, and $D \in \mathbb R^{m \times n}$ is a random direction from the distribution ${\rm Unif}(\mathcal S^{m \times n})$.

In both the one-point and two-point cases, the exact costs $\mathcal J$, as well as the samples $D$ from  ${\rm Unif}(\mathcal S^{m \times n})$ are necessary to compute the gradient estimators.
However, in practice, it is difficult or even impossible to obtain either of them due to different kinds of errors from function approximations, measurement noises, and external disturbances. Two illustrative examples are given as follows.

\begin{example}
\label{Ex1}
Consider the LQR problem described by (\ref{LTI:sys}) and (\ref{cost}). Since the cost (\ref{cost}) is defined over an infinite-horizon, truncation is necessary to compute $\mathcal J(K; x_0)$, i.e., $\mathcal J(K; x_0)$ is approximated by
$\mathcal J_{\delta}(K; x_0) := \sum_{t=0}^{T_\delta - 1} c_t(x_t, u_t)$,
where $T_\delta \ge 1$ is the rollout length.
Given any $T_\delta$, $\mathcal J(K; x_0) > \mathcal J_{\delta}(K; x_0)$.
\end{example}

\begin{example}
Consider a decentralized LQR problem as studied in \cite{li2021distributed, yan2024distributed}. The global cost is given by $\mathcal J(K) := \lim_{t \to \infty} \frac 1 T \sum_{t=1}^T \mathbb E \Big[\frac 1N \sum_{i = 1}^N c_i(t)\Big]$, where agent $i$ only knows its local cost $c_i(t)$, and communicates with its neighbors over a network. As shown in \cite{li2021distributed},
decentralized implementations introduce errors
both in estimating the global costs, and in generating a sample from ${\rm Unif}(\mathcal S^{m \times n})$.
\end{example}

Inspired by the above observations, it is critical to investigate the robustness of derivative-free methods.
We focus on a perturbed derivative-free algorithm as
\begin{equation}
\label{alg:dev_free:dist}
K_{s+1} = K_s - \eta G(K_s) + \eta E_s,
\end{equation}
where $G(K_s)$ is estimated by one-point or two-point settings, and $E_s$ is the perturbation. The perturbed algorithm is summarized in Algorithm 1.

\begin{algorithm}[t]
\caption{Perturbed Derivative-free Methods} 
\label{alg}
\begin{algorithmic}
\STATE {Initialization: initial state $K_0 \in \mathcal K_{st}$, stepsize $\eta > 0$, smoothing radius $r > 0$ and iteration number $T_s$}
\STATE {\textbf{For} $s \in \{0, 1, \dots T_s - 1\}$~\textbf{do}}
\STATE \hspace{0.2cm} Sample $x_0 \in \mathcal D$ 
and $D_s \sim {\rm Unif}(\mathcal S^{m\times n})$
\STATE 
\vspace{-1.2em}
$$
G(K_s) = \begin{cases}
G_r^1(K_s, D_s; x_0) ~{\rm if~one\text{-}point~setting} \\
G_r^2(K_s, D_s; x_0) ~{\rm if~two\text{-}point~setting} \\
\end{cases}
$$
\vspace{-0.8em}
\STATE 
~~~~$K_{s+1} = K_s - \eta G(K_s) + \eta E_s$
\STATE \textbf{Return}~$K_{T_s}$
\end{algorithmic}
\end{algorithm}

\begin{remark}
Policy optimization methods, including the standard gradient descent, natural policy gradient, and Newton gradient algorithms, have been proposed for LQR problems with convergence guarantees \cite{fazel2018global, hu2023toward, mohammadi2021convergence, fatkhullin2021optimizing}. Their robustness has also been analyzed, as perturbations may arise from the inexact gradient computation \cite{pang2021robust, cui2024robust, cui2024small, sontag2022remarks}.
However, to the best of our knowledge, the robustness of derivative-free methods remains unexplored in existing literature. This motivates us to analyze the performance of Algorithm \ref{alg}.
\end{remark}

\begin{remark}
To estimate the policy gradient $\nabla \mathcal J(K)$ by (\ref{one_point}) or (\ref{two_point}), the approximations get better as the smoothing radius $r$ gets smaller, but a smaller $r$ may result in the estimators with larger variance. On the other hand, it is straightforward that both the stepsize $\eta$ and the perturbation $E_s$ significantly affect the convergence of Algorithm \ref{alg}. This paper aims to determine appropriate values for $r$ and $\eta$, and to bound the perturbation $E_s$, to ensure the trajectory of $K_s$ in Algorithm \ref{alg} converges to an $\epsilon$-accurate optimal solution. 
\end{remark}

\section{MAIN RESULTS}

In this section, we first introduce several useful lemmas, and then, analyze the performance of Algorithm \ref{alg}.

\subsection{Supporting Lemmas}

In light of \cite{fazel2018global, hu2023toward}, 
$\mathcal K_{st}$ is a nonconvex set, and thus,  (\ref{reform}) is a nonconvex problem. However, the cost function $\mathcal J(K)$ enjoys the following properties \cite{fazel2018global, malik2020derivative}.

\begin{lemma}
\label{lem:prop:JK}
Given any linear policy $K \in \mathcal K_{st}$, the following statements hold.
\begin{enumerate}
\item There exist positive scalars $\lambda_K$, $\widetilde \lambda_K$ and $\zeta_K$ such that for all $K'$ satisfying $\Vert K' - K\Vert_F \le \zeta_K$, we have
$$|\mathcal J(K') - \mathcal J(K)| \le \lambda_K \Vert K' - K\Vert_F,$$
and 
$$|\mathcal J(K'; x_0) - \mathcal J(K; x_0)| \le \widetilde \lambda_K \Vert K' - K\Vert_F.$$

\item There exist positive scalars $\phi_K$ and  $\beta_K$ such that for all $K'$ satisfying $\Vert K' - K\Vert_F \le \beta_K$, we have
$$\Vert\nabla \mathcal J(K') - \nabla \mathcal J(K)\Vert_F \le \phi_K \Vert K' - K\Vert_F.$$

\item There exists a universal constant $\mu > 0$ such that
$$
\mu(\mathcal J(K) - \mathcal J(K^*)) \le \Vert\nabla \mathcal J(K) \Vert_F^2, 
$$
where $K^*$ is the optimal policy.
\end{enumerate}
\end{lemma}

\begin{remark}
Lemma \ref{lem:prop:JK} indicates that $\mathcal J(K)$ is locally Lipschitz continuous, has a locally Lipschitz continuous gradient, and satisfies the Polyak-Lojasiewicz (PL) condition.
Note that all the scalars in Lemma \ref{lem:prop:JK}
depend on the policy $K$, the initial state $x_0$, as well as the model parameters in (\ref{LTI:sys}) and (\ref{cost}). Their explicit expressions can be found in \cite{fazel2018global, malik2020derivative}.
\end{remark}

Given the initial policy $K_0 \in \mathcal K_{st}$, we define the initial gap to the optimum as $\Delta_0 := \mathcal J(K_0) - \mathcal J(K^*)$, and consider a sublevel set as
$$
\mathcal G^0 := \{K | \mathcal J(K) - \mathcal J(K^*) \le 10 \Delta_0\}.
$$

If $K \notin \mathcal K_{st}$, then $\mathcal J(K)$ is infinite. However, a projection oracle, mapping $K$ back into $\mathcal K_{st}$, is not accessible in Algorithm~\ref{alg}, and thus, it is significant to control the behavior of $K$ during iterations.
With the help of $\mathcal G^0$,  we will show that $\mathcal J(K)$ remains finite with high probability in the next subsection.
Similar to \cite{malik2020derivative, zhao2023global}, we use the constant $10$ in the definition of $\mathcal G^0$. However, we should remark that this choice is not essential, and can be replaced by other constants without affecting the subsequent analysis.

Based on Lemma \ref{lem:prop:JK}, we define
$$
\phi_0 := \sup_{K \in \mathcal  G^0} \phi_K,
\lambda_0 := \sup_{K \in \mathcal  G^0} \lambda_K, ~{\rm and}~
\rho_0 := \inf_{K \in \mathcal  G^0} \rho_K,
$$
where $\rho_K = \min\{\beta_K, \zeta_K\}$.
These quantities transform the local properties of $\mathcal J$ into its global properties over $\mathcal G^0$. 

To facilitate the analysis, we introduce the following uniform bounds over the set $\mathcal G^0$:
\begin{equation}
\label{def:ginf}
g_\infty := \sup_{K \in \mathcal G^0} \Vert G(K) \Vert_F,
\end{equation}
and
\begin{equation}
\label{def:g2}
g_2 := \sup_{K \in \mathcal G^0} \mathbb E \big[\Vert G(K) - \mathbb E[G(K) \mid K]\Vert_F^2 \big].
\end{equation}

Referring to \cite[Corollaries 9 and 10]{malik2020derivative}, we obtain the next lemma.

\begin{lemma}
\label{lem:g}
~
\begin{enumerate}
\item Consider the one-point estimate given by (\ref{one_point}). 
If $r \le \min\{10 \mathcal J(K_0)/\lambda_0, \rho_0\}$, then
$g_2 \le  [20 C_m \mathcal J(K_0)d/r]^2$, and moreover,
$g_\infty \le 20 C_m \mathcal J(K_0)d/r$.

\item Consider the two-point estimate given by (\ref{two_point}). If $r \le \rho_0$, then 
$g_2 \le  d \lambda_0^2$, and moreover,
$g_\infty \le d \lambda_0$.
\end{enumerate}
\end{lemma}

Given a scalar $r > 0$, we define a smoothed version of $\mathcal J(K)$ as
$\mathcal J_r(K) := \mathbb E[\mathcal J(K + rD)]$,
where $D \sim {\rm Unif}(\mathcal S^{m \times n})$, and the expectation is taken with respect to $D$.
In light of \cite{ghadimi2013stochastic, shamir2017optimal}, the next result holds.
\begin{lemma}
\label{lem:smooth:grad}
Consider $G(K)$ computed by (\ref{one_point}) or (\ref{two_point}). Then the following statements hold.
\begin{enumerate}
\item $\nabla \mathcal J_r(K) = \mathbb E[G(K)]$.

\item $\Vert \nabla \mathcal J_r(K) - \nabla \mathcal J(K)\Vert \le \phi_0 r$.
\end{enumerate}	
\end{lemma}

\subsection{Convergence Analysis}

This subsection establishes the convergence of Algorithm \ref{alg}. 

We define the cost error at iteration $s$ as
$$\Delta_s := \mathcal J(K_s) - \mathcal J(K^*),$$ 
and introduce the stopping time
$$\tau := \min\{s \mid \Delta_s > 10 \Delta_0\},$$
which denotes the first iteration that $K_s$ exits the bounded region $\mathcal G^0$.
Let $\mathcal F_s$ be the $\sigma$-field containing all the randomness up to iteration $s$.

Clearly, (\ref{alg:dev_free:dist}) can be viewed as a stochastic gradient descent method with biased gradient estimation. The next lemma characterizes the change in the cost after the one-step update.

\begin{lemma}
\label{lem:onestep}
Consider Algorithm \ref{alg}, where the perturbations $E_s$ are bounded, i.e., $\Vert E_s\Vert_F \le \delta$ for some $\delta > 0$. The following statements hold.
\begin{enumerate}
\item If $r \le \rho_0$
and $\eta \le \min\big\{{1}/{(4\phi_0)}, {\rho_0}/{(g_\infty + \delta)}\big\}$, then
\begin{equation*}
\begin{aligned}
\mathbb E[\Delta_{s+1} \mid \mathcal F_s] \le&
\Big(1 - \frac {\eta \mu}{4}\Big) \Delta_s 
+ \frac {4 \eta (\phi_0 r + \delta)^2}{\mu \theta_0^2} \\
&+ \eta^2 \phi_0 \big(g_2 + 2(\phi_0 r)^2 + \delta^2 \big),
\end{aligned}
\end{equation*}
where  $\theta_0 = \min\{1/(2\phi_0), \rho_0/\lambda_0\}$.
\item Given any $\epsilon > 0$, suppose that $\eta \le \min\big\{{1}/{(4\phi_0)}, \\
{\rho_0}/{(g_\infty + \delta)}\big\}$, 
$$
r \le \min\Big\{\frac {\mu \theta_0}{16\phi_0}\sqrt{\frac {\epsilon}{15}}, \frac {1}{4\phi_0} \sqrt{\frac{\mu\epsilon}{15}}, \rho_0\Big\},$$
and
$$
\delta  \le \min\Big\{\frac {\mu \theta_0}{16}\sqrt{\frac {\epsilon}{15}},  \frac 12\sqrt{\frac{\mu\epsilon}{30}}\Big\}.
$$ 
Then
\begin{equation*}
\begin{aligned}
\mathbb E[\Delta_{s+1} \mid \mathcal F_s] \le \Big(1 - \frac{\eta \mu}{4}\Big) \Delta_s + \eta^2 \phi_0 g_2 + \frac {\eta \mu \epsilon}{120}.
\end{aligned}
\end{equation*}
\end{enumerate}
\end{lemma}

\emph{Proof:}
Since $r \le \rho_0$ and $\eta \le \rho_0/(g_\infty + \delta)$,
we have
$\Vert K_s \pm r D_s - K_s\Vert_F \le \rho_0$, and moreover,
$$\Vert K_{s+1} - K_s\Vert_F \le \eta \big(\Vert G(K_s) \Vert_F + \Vert E_s \Vert_F \big) \le \eta g_\infty + \eta \delta \le \rho_0.$$

Recalling Lemma \ref{lem:prop:JK} yields
\begin{equation*}
\begin{aligned}
&\mathbb E \big[\mathcal J(K_{s+1}) - \mathcal J(K_s) \mid \mathcal F_s \big] \\
&\le \mathbb E \big[\langle \nabla \mathcal J(K_s), K_{s+1} - K_s \rangle_F + \frac{\phi_0}{2}\Vert K_{s+1} - K_s\Vert_F^2 \mid \mathcal F_s \big] \\
&= \mathbb E \big[\langle \nabla \mathcal J(K_s), -\eta G(K_s) + \eta E_s \rangle_F \\
&~~~~+ \frac{\eta^2 \phi_0}{2}\Vert G(K_s) - E_s\Vert_F^2 \mid \mathcal F_s \big].
\end{aligned}
\end{equation*}

Note that $\Vert G(K_s) - E_s\Vert_F^2 \le 2 \Vert G(K_s) \Vert_F^2 + 2\Vert E_s \Vert_F^2$.
By Lemma \ref{lem:smooth:grad}, we obtain
\begin{equation}
\begin{aligned}
\label{pf:ieq1}
&\mathbb E \big[\mathcal J(K_{s+1}) - \mathcal J(K_s) \mid \mathcal F_s \big] \\
& \le - \eta \Vert \nabla \mathcal J(K_s)\Vert_F^2 + \eta (\phi_0 r + \delta) \Vert \nabla \mathcal J(K_s)\Vert_F  \\
&~~~ + \eta^2 \phi_0 \mathbb E \big[\Vert G(K_s)\Vert_F^2 \mid \mathcal F_s \big] + \eta^2 \phi_0 \delta^2.
\end{aligned}
\end{equation}
It is clear that 
\begin{equation*}
\begin{aligned}
\mathbb E \big[\Vert G(K_s)\Vert_F^2 \mid \mathcal F_s \big] &= {\rm var}(G(K_s)) + \Vert \nabla \mathcal J_r(K_s)\Vert_F^2 \\
& \le g_2 + 2\Vert \nabla \mathcal J(K_s)\Vert_F^2 + 2(\phi_0 r)^2.
\end{aligned}
\end{equation*}
Let $\theta$ be a scalar such that $\theta \Vert \nabla \mathcal J(K_s)\Vert_F \le \rho_0$. Then
\begin{equation*}
\begin{aligned}
(\theta - \theta^2 \phi_0/2) \Vert \nabla \mathcal J(K_s)\Vert_F^2 &\le 
\mathcal J(K_s) - \mathcal J(K_s - \theta \nabla \mathcal J(K_s))\\
&\le \mathcal J(K_s) - \mathcal J(K^*).
\end{aligned}
\end{equation*}

Set $\theta = \theta_0 = \min\{1/(2\phi_0), \rho_0/\lambda_0\}$.
Resorting to Lemma \ref{lem:prop:JK} and (\ref{pf:ieq1}), we obtain
\begin{equation*}
\begin{aligned}
&\mathbb E[\Delta_{s+1} - \Delta_s \mid \mathcal F_s] \\
& \le - \eta \Vert \nabla \mathcal J(K_s)\Vert_F^2 + \eta (\phi_0 r + \delta) \frac {2}{\theta_0} \Delta_s^{1/2}  \\
&~~~ + \eta^2 \phi_0 \big(g_2 + 2\Vert \nabla \mathcal J(K_s)\Vert_F^2 + 2(\phi_0 r)^2 + \delta^2 \big).
\end{aligned}
\end{equation*}
Furthermore,
$$\eta (\phi_0 r + \delta) \frac {2}{\theta_0} \Delta_s^{1/2} \le \frac {\eta \mu}{4} \Delta_s + \frac {4 \eta (\phi_0 r + \delta)^2}{\mu \theta_0^2}.$$

Since $\eta \le 1/(4\phi_0)$, it follows from Lemma \ref{lem:prop:JK} iii) that
\begin{equation}
\begin{aligned}
\label{pf:ieq2}
&\mathbb E[\Delta_{s+1} - \Delta_s \mid \mathcal F_s] \\
& \le - \frac {\eta \mu}{4} \Delta_s + \frac {4 \eta (\phi_0 r + \delta)^2}{\mu \theta_0^2}  + \eta^2 \phi_0 \big(g_2 + 2(\phi_0 r)^2 + \delta^2 \big).
\end{aligned}
\end{equation}
Therefore, part i) holds.

Note that $(\phi_0 r + \delta)^2 \le 2(\phi_0 r)^2+ 2 \delta^2$.
Substituting $r$ and $\delta$ of part ii) into (\ref{pf:ieq2}), it can be verified that
\begin{equation*}
\begin{aligned}
\label{pf:J:ine3}
\mathbb E[\Delta_{s+1} - \Delta_s \mid \mathcal F_s]  \le - \frac {\eta \mu}{4} \Delta_s + \eta^2 \phi_0 g_2 + \frac {\eta \mu \epsilon}{120}.
\end{aligned}
\end{equation*}
Thus, the conclusion follows, and this completes the proof.
$\hfill\square$

\begin{remark}
Under biased gradients estimators and perturbations,
Lemma \ref{lem:onestep} i) gives the expected improvement in the cost after one-step iteration,  while ii) further refines the improvement under tighter conditions. The results indicate that $\Delta_s$ decreases geometrically with an error term.
\end{remark}

In the following, we extend the result of Lemma \ref{lem:onestep} ii). Moreover, we show that all the policies $\{K_s\}_{s=1}^{T_s}$  lie in $\mathcal G^0$ with high probability, and thus, belong to $\mathcal K_{st}$ with high probability.

\begin{lemma}
\label{lem:prob}
Consider Algorithm \ref{alg}, where $\Vert E_s\Vert_F \le \delta$ for some $\delta > 0$. 
Let $\epsilon$ be a given error tolerance  such that  $\epsilon \log(120 \Delta_0 / \epsilon) < 12 \Delta_0$.
Suppose that 
$$\eta \le \min\Big\{\frac {\mu \epsilon}{480 \phi_0 g_2}, \frac{1}{4\phi_0}, \frac{\rho_0}{g_\infty + \delta}\Big\},$$
and the parameters $r$ and $\delta$ satisfy the conditions as those in Lemma \ref{lem:onestep} ii).
If $T_s = \frac {4}{\eta\mu}\log \big(\frac {120 \Delta_0}{\epsilon}\big)$,
then
$$\mathbb E[\Delta_{T_s} \cdot 1_{\tau > T_s}] \le {\epsilon}/{20},$$
and furthermore, $\{\tau > T_s\}$ occurs with probability at least $4/5$.
\end{lemma}

\emph{Proof:}
It is straightforward that
$$
\mathbb E[\Delta_{s+1} 1_{\tau > s+1} | \mathcal F_s] \le 
\mathbb E[\Delta_{s+1} 1_{\tau > s} | \mathcal F_s] =
\mathbb E[\Delta_{s+1} | \mathcal F_s] 1_{\tau > s}.
$$

If $\tau > s$, recalling Lemma \ref{lem:onestep} ii) gives 
\begin{equation*}
\begin{aligned}
\mathbb E[\Delta_{s+1} \mid \mathcal F_s] \le \Big(1 - \frac{\eta \mu}{4}\Big)\Delta_s + \eta^2 \phi_0 g_2 + \frac {\eta \mu \epsilon}{120}.
\end{aligned}
\end{equation*}
In the case of $\tau \le s$, we have
$$
\mathbb E[\Delta_{s+1} \mid \mathcal F_s] 1_{\tau > s} = 0.
$$
As a consequence, we obtain
\begin{equation*}
\begin{aligned}
\mathbb E[\Delta_{s+1} \mid \mathcal F_s] 1_{\tau > s} \le \Big(1 - \frac{\eta \mu}{4}\Big)\Delta_s 1_{\tau > s} + \eta^2 \phi_0 g_2 + \frac {\eta \mu \epsilon}{120}.
\end{aligned}
\end{equation*}

Taking expectations over the $\sigma$-field $\mathcal F_s$, we derive
\begin{equation*}
\begin{aligned}
&\mathbb E[\Delta_{s+1} 1_{\tau > s + 1}]  \\
&\le \Big(1 \!-\! \frac{\eta \mu}{4}\Big)^{s+1}\Delta_0 
+ \Big(\eta^2 \phi_0 g_2 + \frac {\eta\mu\epsilon}{120} \Big)\sum_{j=0}^s \Big(1 \!-\! \frac{\eta \mu}{4}\Big)^j \\
&\le \Big(1 - \frac{\eta \mu}{4}\Big)^{s+1}\Delta_0
+ 4 \frac{\eta}{\mu} \phi_0 g_2 + \frac {\epsilon}{30}.
\end{aligned}
\end{equation*}

Let $T_s = \frac {4}{\eta\mu}\log \big(\frac {120 \Delta_0}{\epsilon}\big)$ and $s = T_s -1$. Since $\eta \le {\mu\epsilon}/{(480 \phi_0 g_2)}$, it holds that
$$\mathbb E[\Delta_{T_s} \cdot 1_{\tau > T_s}] \le  {\epsilon}/{20}.$$

Next, we show that $\mathbb P(\tau > T_s) \ge 4/5$ by constructing a super-martingale.

Define 
$$Y_s = \Delta_{\tau \land s} + (2T_s - s)\Big( \eta^2 \phi_0 g_2 + \frac {\eta \mu \epsilon}{120}\Big), ~s \in \{0, 1, \dots, 2T_s\},$$ 
where ${\tau \land s} = \min\{\tau, s\}$. Then
\begin{equation*}
\begin{aligned}
\label{pf:J:ine4}
\mathbb E[Y_{s+1} \mid \mathcal F_s] &= \mathbb E[\Delta_{\tau \land (s+1)} 1_{\tau \le s} \mid \mathcal F_s] \\
& + \mathbb E[\Delta_{\tau \land (s+1)} 1_{\tau > s} \mid \mathcal F_s]	\\
&+\big(2T_s - (s+1)\big)\Big(\eta^2 \phi_0 g_2 + \frac {\eta \mu \epsilon}{120}\Big).
\end{aligned}
\end{equation*}

Note that
$$
\mathbb E \big[\Delta_{\tau \land (s+1)} 1_{\tau \le s} \mid \mathcal F_s \big] 
=\mathbb E \big[\Delta_{\tau \land s} 1_{\tau \le s} \mid \mathcal F_s \big]
= \Delta_{\tau \land s} 1_{\tau \le s},
$$
and moreover,
\begin{equation*}
\begin{aligned}
\mathbb E[\Delta&_{\tau \land (s+1)} 1_{\tau > s} \mid \mathcal F_s] = \mathbb E[\Delta_{s+1} \mid \mathcal F_s] 1_{\tau > s} \\
& \le \Big(1 - \frac{\eta \mu}{4}\Big)\Delta_s 1_{\tau > s} + \Big(\eta^2 \phi_0 g_2 + \frac {\eta \mu \epsilon}{120}\Big) 1_{\tau > s} \\
& \le \Big(1 - \frac{\eta \mu}{4}\Big)\Delta_{\tau \land s} 1_{\tau > s} + \eta^2 \phi_0 g_2 + \frac {\eta \mu \epsilon}{120}.
\end{aligned}
\end{equation*}

Therefore,
\begin{equation*}
\begin{aligned}
\mathbb E[Y&_{s+1} \mid \mathcal F_s] \le
\Delta_{\tau \land s} 1_{\tau \le s} + 
\big(1 - \frac{\eta \mu}{4}\big)\Delta_{\tau \land s} 1_{\tau > s} \\
&+ \eta^2 \phi_0 g_2 + \frac {\eta \epsilon \mu}{120} + (2T_s - (s+1))\Big(\eta^2 \phi_0 g_2 + \frac {\eta\mu \epsilon}{120}\Big) \\
\le& \Delta_{\tau \land s} + (2T_s - s)\Big(\eta^2 \phi_0 g_2 + \frac {\eta \mu \epsilon}{120}\Big) = Y_s,
\end{aligned}
\end{equation*}
and furthermore, $Y_s$ is a super-martingale.

Resorting to the Markov's inequality and $\eta \le {\mu \epsilon}/{(480 \phi_0 g_2)}$, we have
\begin{equation*}
\begin{aligned}
\mathbb P(\max\nolimits&_{s\in \{0, \dots, 2T_s\}} Y_s \ge \nu) \le \frac {\mathbb E[Y_0]}{\nu} \\
&= \frac {1}{\nu}\Big(\Delta_0 + 2T_s\Big(\eta^2 \phi_0 g_2 + \frac {\eta  \mu \epsilon}{120}\Big)\Big)\\
& \le \frac {1}{\nu}\Big(\Delta_0 + \frac{\epsilon}{12}\log(120 \Delta_0/\epsilon)\Big).
\end{aligned}
\end{equation*}

By setting $\nu = 10 \Delta_0$, we conclude that
$\{\tau > T_s\}$ occurs with probability at least $4/5$ if
$\epsilon \log(120 \Delta_0 / \epsilon) < 12 \Delta_0$.
This completes the proof.
$\hfill\square$

We are now ready to establish the convergence of Algorithm \ref{alg} as follows.

\begin{theorem}
\label{thm:conv}
Consider Algorithm \ref{alg}, where $\Vert E_s\Vert_F \le \delta$ for some $\delta > 0$. Let $\epsilon$ be a given error tolerance such that $\epsilon \log(120 \Delta_0 / \epsilon) < 12 \Delta_0$. The following statements hold.
\begin{enumerate}
\item Suppose that the policy gradient is estimated by the one-point setting (\ref{one_point}),
\begin{equation*}
\begin{aligned}
\eta \le \min\Big\{\frac {\mu r^2 \epsilon}{480 \phi_0 \big[20C_m \mathcal J(K_0)d\big]^2}, 
\frac{1}{4\phi_0},& \\
\frac{\rho_0r}{20 C_m \mathcal J(K_0)d + \delta r}\Big\}&,
\end{aligned}
\end{equation*}
\begin{equation*}
\begin{aligned}
r &\le \min\Big\{\frac {\mu \theta_0}{16\phi_0}\sqrt{\frac {\epsilon}{15}}, \frac {1}{4\phi_0} \sqrt{\frac{\mu\epsilon}{15}}, \rho_0, \frac {10 \mathcal J(K_0)}{\lambda_0}\Big\},&
\end{aligned}
\end{equation*}
and
$$
\delta  \le \min\Big\{\frac {\mu \theta_0}{16}\sqrt{\frac {\epsilon}{15}},  \frac 12\sqrt{\frac{\mu\epsilon}{30}}\Big\}.
$$ 
Then $\mathcal J(K_{T_s}) - \mathcal J(K^*) < \epsilon$
with probability at least $3/4$, where $T_s = \frac {4}{\eta\mu}\log \big(\frac {120 \Delta_0}{\epsilon}\big)$.

\item Suppose that the policy gradient is estimated by the two-point setting (\ref{two_point}),
$$\eta \le \min\Big\{\frac {\mu\epsilon}{480 \phi_0 d \lambda_0^2}, \frac{1}{4\phi_0}, \frac{\rho_0}{d \lambda_0 + \delta}\Big\},$$
$$
r \le \min\Big\{\frac {\mu \theta_0}{16\phi_0}\sqrt{\frac {\epsilon}{15}}, \frac {1}{4\phi_0} \sqrt{\frac{\mu\epsilon}{15}}, \rho_0\Big\},$$
and
$$
\delta  \le \min\Big\{\frac {\mu \theta_0}{16}\sqrt{\frac {\epsilon}{15}},  \frac 12\sqrt{\frac{\mu\epsilon}{30}}\Big\}.
$$ 
Then $\mathcal J(K_{T_s}) - \mathcal J(K^*) < \epsilon$
with probability at least $3/4$, where $T_s = \frac {4}{\eta\mu}\log \big(\frac {120 \Delta_0}{\epsilon}\big)$.
\end{enumerate}
\end{theorem}

\emph{Proof:}
With the parameter settings of Lemma \ref{lem:prob}, we have
\begin{equation*}
\begin{aligned}
\mathbb P(\Delta_{T_s} \ge \epsilon) &\le \mathbb P(\Delta_{T_s} \cdot 1_{\tau > T_s} \ge \epsilon) + 
\mathbb P(1_{\tau \le T_s}) \\
& \le \frac 1\epsilon \mathbb E[\Delta_{T_s} \cdot 
1_{\tau > T_s}] + \mathbb P(1_{\tau \le T_s}) \\
&\le 1/20 + 1/5 = 1/4,
\end{aligned}
\end{equation*}
where the second inequality comes from the Markov inequality, and the last inequality follows from Lemma \ref{lem:prob}.
As a consequence, we obtain
$\mathcal J(K_{T_s}) - \mathcal J(K^*) < \epsilon$ with probability at least $3/4$.

For the one-point setting,
the bounds on $g_2$ and $g_\infty$ follow from Lemma \ref{lem:g}~i).
Plugging the bounds into 
$$\eta \le \min\Big\{\frac {\mu \epsilon}{480 \phi_0 g_2}, \frac{1}{4\phi_0}, \frac{\rho_0}{g_\infty + \delta}\Big\},$$
part i) holds.
By a similar procedure for the two-point setting, part ii) follows, and this completes the proof.
$\hfill\square$

\begin{remark}
For any pre-specified error tolerance $\epsilon$, Theorem \ref{thm:conv} indicates that by selecting proper stepsize $\eta$ and smoothing radius $r$, and by bounding the perturbation $E_s$, the trajectory of $K_s$ in Algorithm \ref{alg} can reach an $\epsilon$-accurate optimal solution with high probability. Note that the probability may be further improved by a more refined analysis. 
\end{remark}

\begin{remark}
Under the one-point and two-point settings, Table \ref{tab:summary} summarizes the choices of the stepsize $\eta$ and the smoothing radius $r$, the perturbation bound $\delta$, as well as the sample complexity $T_s$. A faster convergence rate is achieved by two-point gradient estimators due to the reduced variance $g_2$, as shown in Lemma \ref{lem:g}.
\end{remark}

\begin{table}[t]
\centering
\caption{Complexity of derivative-free methods}
\setlength{\tabcolsep}{4pt} 
\label{tab:summary}
\begin{tabular}{ccccc}
\toprule
Parameters & Smoothing radius  & Stepsize & Perturbation  & Complexity  \\
settings & $r$ & $\eta$ & $\delta$ & $T_s$ \\
\midrule
One-point & $\mathcal O(\epsilon^{\frac 12})$ & $\mathcal O(\epsilon^2)$ & $\mathcal O(\epsilon^{\frac 12})$ & 
$\mathcal O \big(\frac {1}{\epsilon^2}\log(\frac 1 \epsilon)\big)$  \\
\rule{0pt}{12pt}
Two-point  &  $\mathcal O(\epsilon^{\frac 12})$ & $\mathcal O(\epsilon)$ & $\mathcal O(\epsilon^{\frac 12})$ & $ \mathcal O\big(\frac 1\epsilon\log(\frac 1 \epsilon) \big)$  \rule{0pt}{12pt} \\ 
\bottomrule
\end{tabular}
\end{table}

As discussed in Example \ref{Ex1}, truncation is necessary to compute the cost function of an infinite-horizon LQR problem and inevitably introduces approximation errors.
The next corollary provides the rollout length for Algorithm \ref{alg} to achieve an $\epsilon$-accurate optimal policy.

\begin{corollary}
Consider solving the LQR problem described by (\ref{LTI:sys}) and (\ref{cost}) using Algorithm \ref{alg}, where the infinite-horizon cost $\mathcal J(K;x_0)$ is approximated by a finite-horizon cost as $\mathcal J_\delta(K;x_0) := \sum_{t = 0}^{T_\delta - 1} c_t(x_t, u_t)$.
Let the parameters $\eta$, $r$, $\epsilon$ and $T_s$ satisfy the requirements as those in Theorem \ref{thm:conv}.
Suppose that the rollout length $T_\delta$ satisfies
\begin{equation*}
\begin{aligned}
T_\delta \ge \frac {1}{2(1 - \gamma)} \max\Big\{
&\log \Big[\frac {16\sqrt{15} M \beta d  \Vert x_0\Vert_2^2} {(1 - \gamma^2) r \mu \theta_0 \sqrt{\epsilon}} \Big], \\
&\log \Big[\frac {2\sqrt{30} M \beta d \Vert x_0\Vert_2^2} {(1 - \gamma^2) r \sqrt{\mu \epsilon}}\Big]
\Big\},
\end{aligned}
\end{equation*}
where $\beta = \max_{K \in \mathcal G^0}\Vert Q + K^\top R K\Vert_2$, $\gamma \in (0, 1)$ and $M > 0$ are constants. 
Then it holds that 
$\mathcal J(K_{T_s}) - \mathcal J(K^*) < \epsilon$
with probability at least $3/4$.
\end{corollary}

\emph{Proof:}
For dynamics (\ref{LTI:sys}), there exist constants 
$\gamma \in (0, 1)$ and $M > 0$ such that for all  $K \in \mathcal G^0$, it holds that 
$$\Vert x_t\Vert_2^2 \le M \gamma^{2t} \Vert x_0\Vert_2^2.$$
Define $\beta := \max_{K \in \mathcal G^0}\Vert Q + K^\top R K\Vert_2$.
Then 
\begin{equation*}
\begin{aligned}
c_t(x_t, u_t) &= \Vert x_t^\top (Q + K^\top R K) x_t \Vert_2 \\
&\le \beta \Vert x_t\Vert^2_2 \le M \beta \gamma^{2t}  \Vert x_0\Vert^2_2.
\end{aligned}
\end{equation*}

Consequently, 
\begin{equation*}
\begin{aligned}
\mathcal J(K) &- \mathcal J_\delta(K) = 
\sum\nolimits_{t = T_{\delta}}^\infty c_t(x_t, u_t) \\
&\le \sum\nolimits_{t = T_{\delta}}^\infty  M \beta \gamma^{2t}\Vert x_0\Vert^2_2
= \frac{M \beta \gamma^{2T_\delta} \Vert x_0\Vert_2^2} {1-\gamma^2}.
\end{aligned}
\end{equation*}
It is straightforward to verify that $[\mathcal J(K) - \mathcal J_\delta(K)]\frac d r \le \delta$ if
\begin{equation}
\label{pf:ieq3}
T_\delta \ge \frac {1}{2(1-\gamma)} 
\log \Big[\frac {M \beta d \Vert x_0\Vert_2^2} {(1 - \gamma^2) r \delta}\Big],
\end{equation}

For Algorithm \ref{alg}, we consider the perturbations $E_s$ arising from the approximation error in evaluating $\mathcal J(K; x_0)$.
Given a policy $K \in \mathcal K_{st}$, 
it is clear that for both the one-point and two-point settings, $\Vert E_s \Vert_F \le \delta$ if $[\mathcal J(K) - \mathcal J_\delta(K)]\frac d r \le \delta$.

With parameter settings as in Theorem \ref{thm:conv}, $\mathcal J(K_{T_s}) - \mathcal J(K^*) < \epsilon$
with probability at least $3/4$. 
Thus, plugging 
\begin{equation*}
\delta  \le \min\Big\{\frac {\mu \theta_0}{16}\sqrt{\frac {\epsilon}{15}},  \frac 12\sqrt{\frac{\mu\epsilon}{30}}\Big\}
\end{equation*}
into (\ref{pf:ieq3}), the conclusion follows. The proof is completed.
$\hfill\square$

\section{NUMERICAL SIMULATIONS}

In this section, we validate the theoretical results through numerical simulations.

Consider the parameters of (\ref{LTI:sys}) and (\ref{cost}) given by
\begin{equation*}
\begin{aligned}
A = \begin{bmatrix}
2  &3 \\
1  &2
\end{bmatrix}, ~~
B = \begin{bmatrix}
1 \\
-1
\end{bmatrix}, ~~
Q = \begin{bmatrix}
2  &-1 \\
-1  &2
\end{bmatrix}, ~~
R = 2
\end{aligned}
\end{equation*}
and the initial state $x_0$ is sampled uniformly at random from the canonical basis vectors.
We solve the problem by Algorithm \ref{alg}.
It can be verified that Assumptions \ref{ass:mat} and \ref{ass:noise} hold.

\begin{figure}[htp]
\centering
\includegraphics[scale=0.33]{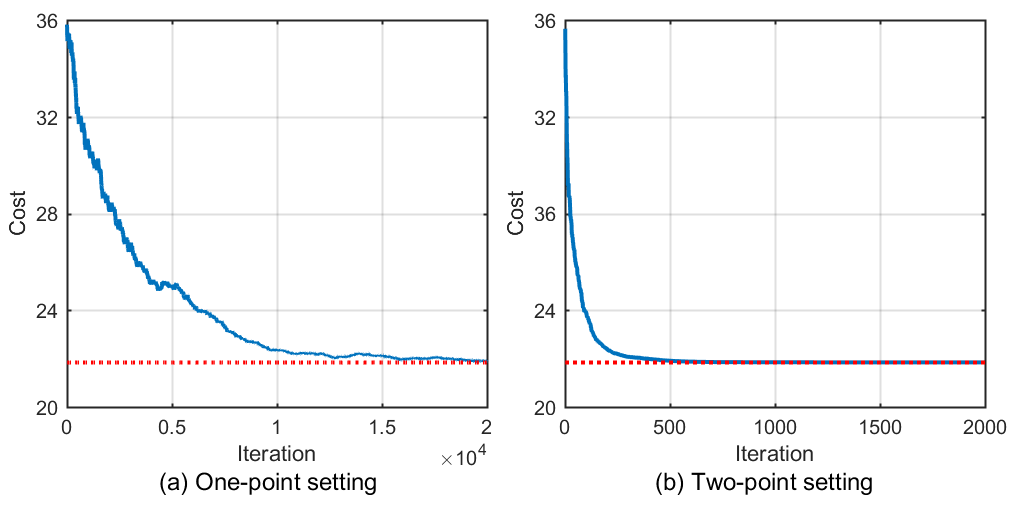}
\caption{Trajectories of $\mathcal J(K_s)$ under both the one-point and two-point settings, where the dashed line shows the cost under the optimal policy. 
(a) One-point setting: $\eta = 2 * 10^{-6}$, $r = 0.1$ and $\delta = 10^{-3}$;
(b) Two-point setting: $\eta = 10^{-4}$, $r = 0.12$ and $\delta = 0.1$.}
\label{Fig1}
\end{figure}

Fig. \ref{Fig1} shows the trajectories of $\mathcal J(K_s)$ under small perturbations. The results indicate that the trajectories of $K_s$ converge to a neighborhood of the optimum.

\begin{figure}[htp]
\centering
\includegraphics[scale=0.33]{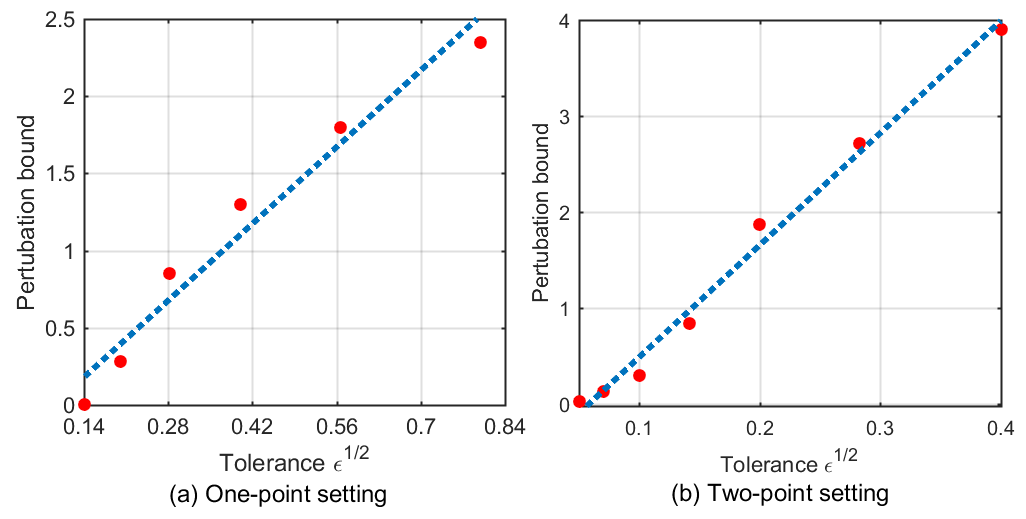}
\caption{The perturbation bound $\delta$ required to ensure Algorithm \ref{alg} achieves an $\epsilon$-accurate optimal policy.
}
\label{Fig2}
\end{figure}

Given the stepsize $\eta$ and the smoothing radius $r$, Fig. \ref{Fig2} presents the perturbation bound required to ensure Algorithm \ref{alg} achieves an $\epsilon$-accurate optimal policy.
The results suggest that the perturbation bound is nearly proportional to $\sqrt{\epsilon}$, which verifies our theoretical findings, as shown in Table \ref{tab:summary}.

\section{CONCLUSION}

In this paper, we analyzed the impact of perturbations on derivative-free optimization methods for infinite-horizon LQR problems. We demonstrated that, under sufficiently small perturbations, these methods can converge to an arbitrarily small neighborhood of the optimal solution through an appropriate selection of algorithmic parameters. We derived explicit bounds on the permissible perturbation magnitude and established the sample complexity required to ensure convergence guarantees. Furthermore, we characterized the necessary rollout length for reliable cost estimation in the derivative-free setting. Finally, we validated our theoretical findings through illustrative numerical simulations.

\section*{References}

\vspace{-1.5em}
\bibliographystyle{IEEEtran}
\bibliography{references,IDS_Publications_03112025}

\end{document}